\documentclass[12pt]{amsart}       
\usepackage{txfonts}
\usepackage{amssymb}
\usepackage{eucal}
\usepackage{graphicx}
\usepackage{amsmath}
\usepackage{amscd}
\usepackage[all]{xy}           

\usepackage{amsfonts,latexsym}
\usepackage{xspace}
\usepackage{epsfig}
\usepackage{float}
\usepackage{color}
\usepackage{fancybox}
\usepackage{colordvi}
\usepackage{multicol}
\usepackage{colordvi}
\usepackage[active]{srcltx} 
\usepackage[colorlinks,final,backref=page,hyperindex,hypertex]{hyperref}

\newcommand{\nc}{\newcommand}
\newcommand{\delete}[1]{}

\nc{\mlabel}[1]{\label{#1}}  
\nc{\mcite}[1]{\cite{#1}}  
\nc{\mref}[1]{\ref{#1}}  
\nc{\mbibitem}[1]{\bibitem{#1}} 

\delete{
\nc{\mlabel}[1]{\label{#1}  
{\hfill \hspace{1cm}{\small\tt{{\ }\hfill(#1)}}}}
\nc{\mcite}[1]{\cite{#1}{\small{\tt{{\ }(#1)}}}}  
\nc{\mref}[1]{\ref{#1}{{\tt{{\ }(#1)}}}}  
\nc{\mbibitem}[1]{\bibitem[\bf #1]{#1}} 
}

\newtheorem{theorem}{Theorem}[section]
\newtheorem{prop}[theorem]{Proposition}

\theoremstyle{definition}
\newtheorem{defn}[theorem]{Definition}
\newtheorem{prop-def}{Proposition-Definition}[section]

\newtheorem{remark}[theorem]{Remark}

\newtheorem{tempex}[theorem]{Example}
\newtheorem{tempexs}[theorem]{Examples}
\newtheorem{temprmk}[theorem]{Remark}
\newtheorem{tempexer}{Exercise}[section]
\newenvironment{exam}{\begin{tempex}\rm}{\end{tempex}}

\nc{\vsa}{\vspace{-.1cm}} \nc{\vsb}{\vspace{-.2cm}}
\nc{\vsc}{\vspace{-.3cm}} \nc{\vsd}{\vspace{-.4cm}}
\nc{\vse}{\vspace{-.5cm}}

\nc{\NS}{U_{NS}}
\nc{\FN}{F_{\mathrm Nij}}
\nc{\dfgen}{V} \nc{\dfrel}{R}
\nc{\dfgenb}{\vec{v}} \nc{\dfrelb}{\vec{r}}
\nc{\dfgene}{v} \nc{\dfrele}{r}
\nc{\dfop}{\odot}
\nc{\dfoa}{\dfop^{(1)}} \nc{\dfob}{\dfop^{(2)}}
\nc{\dfoc}{\dfop^{(3)}} \nc{\dfod}{\dfop^{(4)}}
\nc{\mapm}[1]{\lfloor\!|{#1}|\!\rfloor}
\nc{\cmapm}[1]{\frakC(#1)}
\nc{\red}{\mathrm{Red}}
\nc{\cm}{C}
\nc{\supp}{\mathrm{Supp}}
\nc{\lex}{\mathrm{lex}}

\nc{\disp}[1]{\displaystyle{#1}}
\nc{\bin}[2]{ (_{\stackrel{\scs{#1}}{\scs{#2}}})}  
\nc{\binc}[2]{ \left (\!\! \begin{array}{c} \scs{#1}\\
    \scs{#2} \end{array}\!\! \right )}  
\nc{\bincc}[2]{  \left ( {\scs{#1} \atop
    \vspace{-.5cm}\scs{#2}} \right )}  
\nc{\sarray}[2]{\begin{array}{c}#1 \vspace{.1cm}\\ \hline
    \vspace{-.35cm} \\ #2 \end{array}}
\nc{\bs}{\bar{S}} \nc{\ep}{\epsilon}
\nc{\dbigcup}{\stackrel{\bullet}{\bigcup}}
\nc{\la}{\longrightarrow} \nc{\cprod}{\ast} \nc{\rar}{\rightarrow}
\nc{\dar}{\downarrow} \nc{\labeq}[1]{\stackrel{#1}{=}}
\nc{\dap}[1]{\downarrow \rlap{$\scriptstyle{#1}$}}
\nc{\uap}[1]{\uparrow \rlap{$\scriptstyle{#1}$}}
\nc{\defeq}{\stackrel{\rm def}{=}} \nc{\dis}[1]{\displaystyle{#1}}
\nc{\dotcup}{\ \displaystyle{\bigcup^\bullet}\ }
\nc{\sdotcup}{\tiny{ \displaystyle{\bigcup^\bullet}\ }}
\nc{\fe}{\'{e}}
\nc{\hcm}{\ \hat{,}\ } \nc{\hcirc}{\hat{\circ}}
\nc{\hts}{\hat{\shpr}} \nc{\lts}{\stackrel{\leftarrow}{\shpr}}
\nc{\denshpr}{\den{\shpr}}
\nc{\rts}{\stackrel{\rightarrow}{\shpr}} \nc{\lleft}{[}
\nc{\lright}{]} \nc{\uni}[1]{\tilde{#1}} \nc{\free}[1]{\bar{#1}}
\nc{\freea}[1]{\tilde{#1}} \nc{\freev}[1]{\hat{#1}}
\nc{\dt}[1]{\hat{#1}}
\nc{\wor}[1]{\check{#1}}
\nc{\intg}[1]{F_C(#1)}
\nc{\den}[1]{\check{#1}} \nc{\lrpa}{\wr} \nc{\mprod}{\pm}
\nc{\dprod}{\ast_P} \nc{\curlyl}{\left \{ \begin{array}{c} {} \\
{} \end{array}
    \right .  \!\!\!\!\!\!\!}
\nc{\curlyr}{ \!\!\!\!\!\!\!
    \left . \begin{array}{c} {} \\ {} \end{array}
    \right \} }
\nc{\longmid}{\left | \begin{array}{c} {} \\ {} \end{array}
    \right . \!\!\!\!\!\!\!}
\nc{\lin}{\call} \nc{\ot}{\otimes}
\nc{\ora}[1]{\stackrel{#1}{\rar}}
\nc{\ola}[1]{\stackrel{#1}{\la}}
\nc{\scs}[1]{\scriptstyle{#1}} \nc{\mrm}[1]{{\rm #1}}
\nc{\margin}[1]{\marginpar{\rm #1}}   
\nc{\dirlim}{\displaystyle{\lim_{\longrightarrow}}\,}
\nc{\invlim}{\displaystyle{\lim_{\longleftarrow}}\,}
\nc{\mvp}{\vspace{0.5cm}}
\nc{\mult}{m}       
\nc{\svp}{\vspace{2cm}} \nc{\vp}{\vspace{8cm}}
\nc{\proofbegin}{\noindent{\bf Proof: }}
\nc{\proofend}{$\blacksquare$ \vspace{0.5cm}}
\nc{\sha}{{\mbox{\cyr X}}}  
\nc{\ncsha}{{\mbox{\cyr X}^{\mathrm NC}}}
\newfont{\scyr}{wncyr10 scaled 550}
\nc{\ssha}{\mbox{\bf \scyr X}}
\nc{\ncshao}{{\mbox{\cyr
X}^{\mathrm NC,\,0}}}
\nc{\shpr}{\diamond}    
\nc{\shprc}{\shpr_c}
\nc{\shpro}{\diamond^0}    
\nc{\shpru}{\check{\diamond}} \nc{\spr}{\cdot}
\nc{\catpr}{\diamond_l} \nc{\rcatpr}{\diamond_r}
\nc{\lapr}{\diamond_a} \nc{\lepr}{\diamond_e} \nc{\sprod}{\bullet}
\nc{\un}{u}                 
\nc{\vep}{\varepsilon} \nc{\labs}{\mid\!} \nc{\rabs}{\!\mid}
\nc{\hsha}{\widehat{\sha}} \nc{\lsha}{\stackrel{\leftarrow}{\sha}}
\nc{\rsha}{\stackrel{\rightarrow}{\sha}} \nc{\lc}{\lfloor}
\nc{\rc}{\rfloor} \nc{\sqmon}[1]{\langle #1\rangle}
\nc{\altx}{\Lambda} \nc{\vecT}{\vec{T}} \nc{\piword}{{\mathfrak P}}
\nc{\lbar}[1]{\overline{#1}}
\nc{\dep}{\mathrm{dep}}

\nc{\mmbox}[1]{\mbox{\ #1\ }}
\nc{\ayb}{\mrm{AYB}} \nc{\mayb}{\mrm{mAYB}} \nc{\cyb}{\mrm{cyb}}
\nc{\ann}{\mrm{ann}} \nc{\Aut}{\mrm{Aut}} \nc{\cabqr}{\mrm{CABQR
}} \nc{\can}{\mrm{can}} \nc{\colim}{\mrm{colim}}
\nc{\Cont}{\mrm{Cont}} \nc{\rchar}{\mrm{char}}
\nc{\cok}{\mrm{coker}} \nc{\dtf}{{R-{\rm tf}}} \nc{\dtor}{{R-{\rm
tor}}}

\nc{\Div}{{\mrm Div}} \nc{\End}{\mrm{End}} \nc{\Ext}{\mrm{Ext}}
\nc{\FG}{\mrm{FG}} \nc{\Fil}{\mrm{Fil}} \nc{\Frob}{\mrm{Frob}}
\nc{\Gal}{\mrm{Gal}} \nc{\GL}{\mrm{GL}} \nc{\Hom}{\mrm{Hom}}
\nc{\hsr}{\mrm{H}} \nc{\hpol}{\mrm{HP}} \nc{\id}{\mrm{id}} \nc{\Id}{\mathrm{Id}}
\nc{\im}{\mrm{im}} \nc{\incl}{\mrm{incl}} \nc{\Loday}{\mrm{ABQR}\
} \nc{\length}{\mrm{length}} \nc{\LR}{\mrm{LR}} \nc{\mchar}{\rm
char} \nc{\pmchar}{\partial\mchar} \nc{\map}{\mrm{Map}}
\nc{\MS}{\mrm{MS}} \nc{\OS}{\mrm{OS}} \nc{\NC}{\mrm{NC}}
\nc{\rba}{\rm{Rota-Baxter algebra}\xspace}
\nc{\rbas}{\rm{Rota-Baxter algebras}\xspace}
\nc{\rbw}{\rm{RBW}\xspace}
\nc{\rbws}{\rm{RBWs}\xspace}
\nc{\rbadj}{\rm{RB}\xspace}
\nc{\mpart}{\mrm{part}} \nc{\ql}{{\QQ_\ell}} \nc{\qp}{{\QQ_p}}
\nc{\rank}{\mrm{rank}} \nc{\rcot}{\mrm{cot}} \nc{\rdef}{\mrm{def}}
\nc{\rdiv}{{\rm div}} \nc{\rtf}{{\rm tf}} \nc{\rtor}{{\rm tor}}
\nc{\res}{\mrm{res}} \nc{\SL}{\mrm{SL}} \nc{\Spec}{\mrm{Spec}}
\nc{\tor}{\mrm{tor}} \nc{\Tr}{\mrm{Tr}}
\nc{\mtr}{\mrm{tr}}

\nc{\ab}{\mathbf{Ab}} \nc{\Alg}{\mathbf{Alg}}
\nc{\Bax}{\mathbf{CRB}} \nc{\Algo}{\mathbf{Alg}^0}
\nc{\cRB}{\mathbf{CRB}} \nc{\cRBo}{\mathbf{CRB}^0}
\nc{\RBo}{\mathbf{RB}^0} \nc{\BRB}{\mathbf{RB}}
\nc{\Dend}{\mathbf{DD}} \nc{\bfk}{{\bf k}} \nc{\bfone}{{\bf 1}}
\nc{\base}[1]{{a_{#1}}} \nc{\Cat}{\mathbf{Cat}}
 \nc{\DN}{\mathbf{DN}}
\nc{\NA}{\mathbf{NA}}
\nc{\SDN}{\mathbf{SDN}}
\nc{\Diff}{\mathbf{Diff}} \nc{\gap}{\marginpar{\bf
Incomplete}\noindent{\bf Incomplete!!}
    \svp}
\nc{\FMod}{\mathbf{FMod}} \nc{\Int}{\mathbf{Int}}
\nc{\Mon}{\mathbf{Mon}}
\nc{\RB}{\mathbf{RB}} \nc{\remarks}{\noindent{\bf Remarks: }}
\nc{\Rep}{\mathbf{Rep}} \nc{\Rings}{\mathbf{Rings}}
\nc{\Sets}{\mathbf{Sets}} \nc{\bfx}{\mathbf{x}}
\nc{\BA}{{\Bbb A}} \nc{\CC}{{\Bbb C}} 
\nc{\EE}{{\Bbb E}} \nc{\FF}{{\Bbb F}} \nc{\GG}{{\Bbb G}}
\nc{\HH}{{\Bbb H}} \nc{\LL}{{\Bbb L}} \nc{\NN}{{\Bbb N}}
\nc{\QQ}{{\Bbb Q}} \nc{\RR}{{\Bbb R}} \nc{\TT}{{\Bbb T}}
\nc{\VV}{{\Bbb V}} \nc{\ZZ}{{\Bbb Z}}

\nc{\cala}{{\mathcal A}} \nc{\calb}{{\mathcal B}}
\nc{\calc}{{\mathcal C}}
\nc{\cald}{{\mathcal D}} \nc{\cale}{{\mathcal E}}
\nc{\calf}{{\mathcal F}} \nc{\calg}{{\mathcal G}}
\nc{\calh}{{\mathcal H}} \nc{\cali}{{\mathcal I}}
\nc{\calj}{{\mathcal J}} \nc{\call}{{\mathcal L}}
\nc{\calm}{{\mathcal M}} \nc{\caln}{{\mathcal N}}
\nc{\calo}{{\mathcal O}} \nc{\calp}{{\mathcal P}}
\nc{\calr}{{\mathcal R}} \nc{\cals}{{\mathcal S}} \nc{\calt}{{\mathcal T}}
\nc{\calw}{{\mathcal W}} \nc{\calx}{{\mathcal X}}
\nc{\CA}{\mathcal{A}}

\nc{\frakA}{{\mathfrak A}}
\nc{\fraka}{{\mathfrak a}}
\nc{\frakB}{{\mathfrak B}}
\nc{\frakb}{{\mathfrak b}}
\nc{\frakC}{{\mathfrak C}}
\nc{\frakd}{{\mathfrak d}}
\nc{\frakF}{{\mathfrak F}}
\nc{\frakg}{{\mathfrak g}}
\nc{\frakm}{{\mathfrak m}}
\nc{\frakM}{{\mathfrak M}}
\nc{\frakMo}{{\mathfrak M}^0}
\nc{\frakP}{{\mathfrak P}}
\nc{\frakp}{{\mathfrak p}}
\nc{\frakS}{{\mathfrak S}}
\nc{\frakSo}{{\mathfrak S}^0}
\nc{\fraks}{{\mathfrak s}}
\nc{\os}{\overline{\fraks}}
\nc{\frakT}{{\mathfrak T}}
\nc{\frakTo}{{\mathfrak T}^0}
\nc{\oT}{\overline{T}}
\nc{\frakX}{{\mathfrak X}}
\nc{\frakXo}{{\mathfrak X}^0}
\nc{\frakx}{{\mathbf x}}
\nc{\frakTx}{\frakT}      
\nc{\frakTa}{\frakT^a}        
\nc{\frakTxo}{\frakTx^0}   
\nc{\caltao}{\calt^{a,0}}   
\nc{\ox}{\overline{\frakx}} \nc{\fraky}{{\mathfrak y}}
\nc{\frakz}{{\mathfrak z}} \nc{\oX}{\overline{X}} \font\cyr=wncyr10

\newcommand{\q}{Q}        
\newcommand{\D}{d}        
\newcommand{\pp}{\Pi}            
\newcommand{\DD}{D}             
\newcommand{\E}{E}                   
\newcommand{\s}{S}                   
\newcommand{\T}{T}                   
\newcommand{\ct}{T}   
 \nc{\ID}{\mathbf{ID}}   
\nc{\CID}{\mathbf{CID}} 
\nc{\fid}[1]{\mathrm{ID}(#1)} 
\nc{\efid}[1]{\mathrm{ID}(#1)^*} 
\nc{\ee}{\vep} \nc{\diffa}[1]{\{#1\}} \nc{\diffs}[1]{\Delta{#1}}

\nc{\tred}[1]{\textcolor{red}{#1}}
\nc{\tgreen}[1]{\textcolor{green}{#1}}
\nc{\tblue}[1]{\textcolor{blue}{#1}} \nc{\li}[1]{\tred{Li:#1 }}
\nc{\xing}[1]{\tblue{Xing:#1 }} \nc{\referee}[1]{\tgreen{Referee:#1 }}

\begin{document}

\title{Constructions of free commutative integro-differential algebras}

\author{Xing Gao}
\address{Department of Mathematics,
    Lanzhou University,
    Lanzhou, Gansu 730000, China}
\email{gaoxing@lzu.edu.cn}
\author{Li Guo}
\address{
Department of Mathematics and Computer Science,
Rutgers University,
Newark, NJ 07102, USA}
\email{liguo@rutgers.edu}

\date{\today}

\begin{abstract}
In this survey, we outline two recent constructions of free commutative integro-differential algebras. They are based on the construction of free commutative Rota-Baxter algebras by mixable shuffles. The first is by evaluations. The second is by the method of Gr\"obner-Shirshov bases.
\end{abstract}

\subjclass[2010]{
16S15, 
13P10, 
16W99, 
12H05, 
47G20 
.}

\keywords{Differential algebra, Rota-Baxter algebra, integro-differential algebra, Gr\"obner-Shirshov basis, free algebra, shuffle product, mixable shuffle product.}

\maketitle

\tableofcontents

\allowdisplaybreaks

\setcounter{section}{0}

\section{Introduction}
In this survey article, we give an outline of the recent constructions of free commutative integro-differential algebras.

The main axiom of integro-differential algebra can be regarded as an algebraic abstraction of the integral by parts formula which involves both derivation and integration. Thus to understand this abstraction better, we first review the abstraction for derivation and for integration.

In this paper, by an algebra we mean a {\em commutative} associative algebra over some commutative ring, unless otherwise specified. A differential algebra is an algebra $R$ together with a linear operator $d\colon R\to R$ that satisfies the following axiom distilled from the Leibniz rule for derivations
$$ d(xy) =d(x)y+xd(y), \text{ for all } x, y\in R.$$
The study of differential algebra began with Ritt's classic work~\cite{Ri1,Ri2}. After the fundamental work of Kolchin~\cite{Ko}, differential algebra has evolved into a vast area of mathematics that is important in both theory~\cite{CGKS,SP} and applications: for instance, in mechanic theorem proving by W.-T. Wu~\cite{Wu,Wu2}. Free (commutative) differential algebras, in the form of differential polynomial algebras (Theorem~\ref{thm:diff}), are essential for studying differential equations, as polynomial algebras are for commutative algebras.

The algebraic study of integrals came much later. In fact the development did not start from an algebraic abstraction of integrals, but from the effort of G. Baxter~\cite{Ba} in 1960 to understand a formula in probability theory. As a result, the concept is not called an integral algebra, but called a (Rota-)Baxter algebra (Eq.~(\mref{eq:rb})) which is the integral counterpart of the derivation, the difference operator, and divided differences (see Eq.~(\ref{eq:diff})). Soon afterwards Rota noticed its importance in combinatorics and promoted its study through research and survey articles (see e.g. ~\cite{Ro1,Ro2}). Independently, Rota-Baxter operators on Lie algebras were found to be closely related to the classical Yang-Baxter equation~\cite{STS}. Since the turn of this century, the theory of Rota-Baxter algebra has experienced rapid development with broad applications in mathematics and physics~\mcite{Bai,Gub,GZ,Ro1,Ro2,STS},  especially noteworthy in the Hopf algebra approach of Connes-Kreimer to renormalization of quantum field theory~\cite{CK,EGK,GZ}. Here again a fundamental role is played by free (commutative) Rota-Baxter algebras that were first constructed by Rota~\cite{Ro1} and Cartier~\cite{Ca}, and then by Guo-Keigher~\cite{GK1} in terms of mixable shuffles (Theorem~\ref{thm:shua}).

The fusion of differential and Rota-Baxter algebras, motivated by algebraic study of calculus as a whole, appeared about five years ago. It is amazing that two structures for this purpose were introduced at about the same time. One is a relatively simple coupling of differential algebra and Rota-Baxter algebra through section axiom (Eq.~(\mref{eq:fft})) that reflects the First Fundamental Theorem of Calculus. It is called differential Rota-Baxter algebra~\cite{GK3}. The other one is a more faithful abstraction of the integration-by-parts formula (see Eq.~(\ref{eq:ibp1})), giving rise to the concept of an integro-differential algebra~\cite{RR} which has generated much interest~\cite{ACPRR,ACPaRR,RRTB,RRTB1}. As suggested in  previous cases, free objects for these algebraic structures are important in their studies. Because of the relative independence of the differential and integral (Rota-Baxter) structures in a differential Rota-Baxter algebra, the free object was constructed by a clear combination of the free objects on the differential and Rota-Baxter sides and were obtained at the same time when the concept was introduced. In contrast, the construction of free integro-differential algebras took longer to achieve. Nevertheless, there are two recent constructions~\cite{GGZ,GRR} and it is the purpose of this paper to give the preliminary background and a summary of these constructions.

Both are based on the construction of free Rota-Baxter algebras by mixable shuffles. Straight from the definition, a free commutative integro-differential algebra can be obtained as the quotient of a free Rota-Baxter algebra modulo the integral-by-parts axiom. By an explicit construction of a free integro-differential algebra, we mean identifying a specific vector space basis of this quotient. Thus we give two such bases in this paper.

After a preliminary Section~2 on the concepts of differential, Rota-Baxter and differential Rota-Baxter algebras as well as operated algebras, and the constructions of their respective free objects, we give the first construction~\cite{GRR} of free integro-differential algebras in Section~3. This construction applies to regular differential algebras, a concept which we also review in Section~3. Common examples of regular differential algebras include differential polynomial algebras and rational functions. The second construction~\cite{GGZ} is given in Section~4. The construction applies the general method of Gr\"obner-Shirshov bases, of which Gr\"obner bases in commutative algebra are special cases, but which apply to many other algebraic structures. We give some details of the method in the case of integro-differential algebras, where we use the ambient algebraic structure of a free Rota-Baxter algebra to establish the Composition-Diamond Lemma.

\section{Definitions and preliminary constructions}

\mlabel{sec:prel} We recall the definitions of algebras with various
differential and integral operators and the constructions of the
free objects in the corresponding categories. Free commutative integro-differential algebras, which are the focus of this survey, will be discussed in later sections.

\subsection{The definitions}
We recall the algebraic structures considered in this paper. We also
introduce variations with nilpotent derivation that will be needed later. Algebras considered in this paper are assumed to be unitary (and commutative), unless explicitly designated as non-unitary.

\begin{defn}
{\rm Let $\bfk$ be a unitary commutative ring. Let $\lambda\in \bfk$
be fixed.
\begin{enumerate}
\item A {\bf differential $\bfk$-algebra of weight $\lambda$} (also
    called a {\bf $\lambda$-differential $\bfk$-algebra}) is an associative $\bfk$-algebra $R$ together with a linear operator
    $\D\colon R\to R$ such that
    \begin{equation}
      \D(xy)=\D(x) y+x \D(y)+ \lambda \D(x)\D(y)\quad \text{for all } x,y\in R, \mlabel{eq:diff}
    \end{equation}
    and
    \begin{equation}
      \D(1)=0.
      \mlabel{eq:diffc}
    \end{equation}
    Such an operator is called a {\bf derivation of weight $\lambda$} or a
    {\bf $\lambda$-derivation}.
\item A {\bf Rota-Baxter $\bfk$-algebra of weight $\lambda$} is an
    associative $\bfk$-algebra $R$ together with a linear operator
    $P\colon R\to R$ such that
\begin{equation}
P(u)P(v)=P(uP(v))+P(P(u)v)+\lambda P(uv) \quad \text{ for all } u, v\in
R.\mlabel{eq:rb}
\end{equation}
\item A {\bf differential Rota-Baxter k-algebra of weight $\lambda$} (also called a {\bf $\lambda$-differential
      Rota-Baxter $\bfk$-algebra}) is a differential $\bfk$-algebra
    $(R,d)$ of weight $\lambda$ with a Rota-Baxter operator $P$ of
    weight $\lambda$ such that
\begin{equation} d\circ P=\id.
\mlabel{eq:fft}
\end{equation}
\item An {\bf integro-differential $\bfk$-algebra of weight
      $\lambda$} (also called a {\bf $\lambda$-integro-differential
      $\bfk$-algebra}) is a differential $\bfk$-algebra $(R,\DD)$ of
    weight $\lambda$ with a linear operator $\pp\colon R \to R$ such
    that
    \begin{equation}
      \DD \circ \pp = \id_R \mlabel{eq:idcomp}
    \end{equation}
    and
    \begin{equation}
      \pp(\DD(x)) \pp(\DD(y)) = \pp(\DD(x)) y + x \pp(\DD(y)) -
      \pp(\DD(xy))  \quad \text{for all } x,y\in R.
      \mlabel{eq:diffbaxter}
    \end{equation}
\end{enumerate}
}
\end{defn}

Eqs.~($\mref{eq:rb}$), ($\mref{eq:idcomp}$) and ($\mref{eq:diffbaxter}$) are called the {\bf Rota-Baxter axiom}, {\bf section axiom} and {\bf hybrid Rota-Baxter axiom}, respectively. It is proved in~\mcite{GRR} that a differential $\bfk$-algebra $(R,\DD)$ with a linear operator $\pp:R\to R$ is an integro-differential algebra if and only if Eq.~(\mref{eq:idcomp}) and the following {\bf integration by parts axioms} hold:
\begin{equation}
x\pp(y) = \pp(\DD(x) \pp(y) ) + \pp(xy) + \lambda \pp(\DD(x) y )
\mlabel{eq:ibp1}
\end{equation}
and
\begin{equation}
\pp(x)y =  \pp(\pp(x) \DD(y) ) + \pp(xy) + \lambda \pp(x
\DD(y) )   \quad \text{for all } x,y\in R.
\mlabel{eq:ibp2}
\end{equation}
These two equations can be regarded as the weighed and noncommutative versions of the classical integration by parts formula in analysis.

\begin{exam}
Let $R = C^{\infty}(\mathbb{R})$
\begin{enumerate}
\item
Fix a $\lambda\in \mathbb{R}$. Define
$$ D_\lambda: R\longrightarrow R, \quad f(x)\mapsto \frac{f(x+\lambda)-f(x)}{\lambda}. $$
Then $D_\lambda$ is a differential operator of weight $\lambda$.
\item
For fixed $a\in \mathbb{R}$, the integral operator
$$ \Pi: R\longrightarrow R, \quad f(x)\mapsto \int_a^x f(t)dt$$
is a Rota-Baxter operator of weight zero.
\item
$D$ be the usual derivation on $R$ and $\Pi$ be the above integral operator.
Then $(R, D, \Pi)$ is a differential Rota-Baxter algebra and an integro-differential algebra of weight 0.
\end{enumerate}
\end{exam}
See~\mcite{Gub,GK3,GRR,Ro2} for more examples.

\subsection{Free differential Rota-Baxter algebras}

We first recall the construction of free commutative differential
algebras and introduce their order $n$ variations. For a set $Y$,
let $C(Y)$ denote the free commutative monoid on $Y$. Thus elements
in $\cm(Y)$ are commutative words from the
alphabet set $Y$, plus the identity $1$. Let $ \bfk [Y]$ be the commutative polynomial algebra generated by $Y$.


\begin{theorem} (\mcite{GGZ,GK3})
\begin{enumerate}
\item
Let $Y$ be a set with a map $d_0\colon Y\to \bfk[Y]$. Extend $d_0$ to $d
\colon \bfk[Y]  \to \bfk[Y]$ as follows. Let $w=u_1\cdots u_k$, where $u_i\in Y$ for $1\leq i\leq k$, be a commutative word from the
alphabet set $Y$. Recursively define
\begin{equation}
      d(w)=d_0(u_1)u_2 \cdots u_k + u_1 d (u_2 \cdots u_k) + \lambda
      d_0(u_1)d(u_2 \cdots u_k).
      \mlabel{eq:prodind}
\end{equation}
Explicitly,
$$
d(w)=\sum_{\emptyset \neq I\subseteq \{1,\cdots,k\}} \lambda ^{|I|-1} \check{d}_I(u_1)\cdots \check{d}_I(u_k), \text{ where }
\check{d}_I(u_i)=\left \{\begin{array}{ll} d_0(u_i), & i\in I, \\ u_i, & i\not\in I. \end{array} \right .
$$
    Further define $d(1)=0$ and then extend $d$ to
    $\bfk[Y]$ by linearity. Then $(\bfk[Y],d)$ is a differential algebra of weight $\lambda$.
\mlabel{it:diffs}
\item
Let $X$ be a set. Let $ Y:=\Delta X: = \{ x^{(n)} \mid x\in X, n\geq
0\}$ with the map $d_0\colon\Delta X\to \Delta X, x^{(n)}\mapsto
x^{(n+1)}$. Then with the extension $d_X:=d$ of $d_0$ as in
Eq.~(\mref{eq:prodind}), $(\bfk\{X\}, d_X) := (\bfk[\Delta X], d_X)$ is the free
commutative differential algebra of weight $\lambda$ on the set $X$.
\mlabel{it:commfreediff}
\item
For a given $n\geq 1$, let $\Delta
X^{(n+1)}:=\left\{x^{(k)}\,\big|\, x\in X, k\geq n+1\right\}$. Then
$\bfk\{X\}\Delta X^{(n+1)}$ is the differential ideal $I_n$ of
$\bfk\{X\}$ generated by the set $\{ x^{(n+1)}\,|\,x\in X\}$.
The quotient differential algebra $\bfk\{X\}/I_n$
has a canonical basis given by $\Delta_n X:=\{x^{(k)}\,|\, k\leq
n\}$, thus giving a differential algebra isomorphism
$\bfk\{X\}/I_n\cong \bfk[\Delta_n X]$ where the differential
structure on the later algebra is given by $d$ in Eq.~(\mref{eq:prodind}), where
\begin{equation}
\mlabel{eq:diffn} d_0(x^{(i)}) = \left\{ \begin{array}{ll} x^{(i+1)},
& 1\leq i\leq n-1, \\ 0, & i=n. \end{array}\right .
\end{equation}
\mlabel{it:diffordn}
\end{enumerate}
\mlabel{thm:diff}
\end{theorem}

We note that in $\bfk[\Delta_n X]$, $d_0^{n+1}(u)=0$ only holds for the variables $x\in X$. For example, when $n=1$, we have $d_0^2(x^2)=2 x^{(1)}\neq 0$.

We next recall the construction of free commutative Rota-Baxter
algebras in terms of mixable shuffles~\mcite{GK1,GK2}. The mixable
shuffle product is shown to be the same as the quasi-shuffle product
of Hoffman~\mcite{EG,GZ,Ho}.  Let $A$ be a commutative
$\bfk$-algebra. Define

\begin{equation}
  \sha (A)= \bigoplus_{k\geq 0} A^{\otimes (k+1)} = A\oplus A^{\otimes
    2}\oplus \cdots.
    \mlabel{eq:freerb}
\end{equation}
Let $\fraka =a_0\ot \cdots \ot a_m\in A^{\ot (m+1)}$ and
$\frakb=b_0\ot \cdots \ot b_n\in A^{\ot (n+1)}$. If $m=0$ or $n=0$,
define
\begin{equation}
  \fraka \shpr \frakb =\left \{\begin{array}{ll}
      (a_0b_0)\ot b_1\ot \cdots \ot b_n, & m=0, n>0,\\
      (a_0b_0)\ot a_1\ot \cdots \ot a_m, & m>0, n=0,\\
      a_0b_0, & m=n=0.
    \end{array} \right .
\end{equation}
If $m>0$ and $n>0$, inductively (on $m+n$) define
\begin{eqnarray}
  \fraka \shpr \frakb & = &
  (a_0b_0)\ot \Big(
  (a_1\ot a_2\ot \cdots \ot a_m) \shpr (1\ot b_1\ot \cdots \ot b_n) \notag \\
  &&
  \qquad \qquad +
  \; (1\ot a_1\ot \cdots \ot a_m) \shpr (b_1\ot \cdots \ot b_n) \mlabel{eq:shpr}\\
  && \qquad \qquad +
  \lambda\, (a_1\ot \cdots \ot a_m) \shpr (b_1\ot \cdots \ot b_n)\Big).
  \notag
\end{eqnarray}
Extending by additivity, we obtain a $\bfk$-bilinear map
\begin{equation*}
  \shpr\colon \sha (A) \times \sha (A) \rar \sha (A).
\end{equation*}
Alternatively,
\begin{equation*}
  \fraka\shpr \frakb=(a_0b_0)\otimes (\lbar{\fraka}
  \ssha_\lambda
  \lbar{\frakb}),
\end{equation*}
where~$\bar{\fraka} = a_1 \otimes \cdots \otimes a_m$, $\bar{\frakb}
= b_1 \otimes \cdots \otimes b_n$ and~$\ssha_\lambda$ is the mixable
shuffle (quasi-shuffle) product of weight
$\lambda$~\mcite{Gub,GK1,Ho}, which specializes to the shuffle
product ${\tiny \ssha}$ when $\lambda=0$.

Define a $\bfk$-linear endomorphism $P_A$ on $\sha (A)$ by assigning
\[ P_A( x_0\otimes x_1\otimes \cdots \otimes x_n) =1_A\otimes
x_0\otimes x_1\otimes \cdots\otimes x_n, \] for all $x_0\otimes
x_1\otimes \cdots\otimes x_n\in A^{\otimes (n+1)}$ and extending by
additivity.  Let $j_A\colon A\rar \sha (A)$ be the canonical
inclusion map.

\begin{theorem} $($\cite{GK1,GK2}$)$ Let $A$ be a commutative $\bfk$-algebra.
\begin{enumerate}
\item
The pair $(\sha(A),P_A)$, together with the
  natural embedding $j_A\colon A\hookrightarrow \sha (A)$, is the free
  commutative Rota-Baxter $\bfk$-algebra on $A$ of weight $\lambda$.
  In other words, for any commutative Rota-Baxter $\bfk$-algebra $(R,P)$ and any
  $\bfk$-algebra map $\varphi\colon A\rar R$, there exists a unique
  Rota-Baxter $\bfk$-algebra homomorphism $\tilde{\varphi}\colon (\sha
  (A),P_A)\rar (R,P)$ such that $\varphi = \tilde{\varphi} \circ
  j_A$ as $\bfk$-algebra homomorphisms.
\item Let $Y$ be a set and let $\bfk[Y]$ be the free
    commutative algebra on $Y$. The pair $(\sha(Y),P_Y):=(\sha (\bfk[Y]),
    P_{\bfk[Y]})$, together with the  natural embedding $j_Y\colon Y\hookrightarrow \bfk[Y]\rightarrow \sha (\bfk[Y])$, is the free
  commutative Rota-Baxter $\bfk$-algebra of weight $\lambda$ on $Y$.
\end{enumerate}
  \mlabel{thm:shua}
\end{theorem}

Since $\shpr$ is compatible with the multiplication in $A$, we will
often suppress the symbol $\shpr$ and simply write $x y$ for
$x\shpr y$ in $\sha (A)$, unless there is a danger of confusion.

A linear basis of $\sha(\bfk[Y])$ is given by
\begin{equation}
\calb(Y):= \left\{ x_0\ot \cdots \ot x_k\,\big|\, x_i\in C(Y),
0\leq i\leq k, k\geq 0\right\}, \mlabel{eq:rbm}
\end{equation}
called the set of {\bf Rota-Baxter monomials} in $Y$. The integer
$\dep(x_0\ot \cdots \ot x_k):=k+1$ is called the {\bf depth} of
$x_0\ot \cdots \ot x_k$. As a convenience, we also write $P$
for $P_{\bfk[Y]}$. Then $1\ot u$ and $P(u)$ stand for the same
element, and we will use both notations synonymously in this paper.

We now put the differential and Rota-Baxter algebra structures
together. Let $(A, d)$ be a commutative differential
$\bfk$-algebra of weight $\lambda$. Extend $d$ to $\sha(A)$ by
\begin{eqnarray*}
  \lefteqn{ d_A(u_0\otimes u_1\otimes\cdots\otimes u_k)}\\
  &=&
  d(u_0)\otimes u_1\otimes \cdots \otimes u_k + u_0u_1\otimes
  u_2 \otimes \cdots \otimes u_k +\lambda d(u_0) u_1\otimes u_2\otimes
  \cdots \otimes u_k, \quad k\geq 0.
\end{eqnarray*}
Note that $d_A$ does not satisfy the Leibniz rule with respect to the tensor product since here a tensor factor means an application of the Rota-Baxter operator $P_A$: $u_0\ot u_1= u_0P_A(u_1)$. Thus
$$ d_A(u_0\ot u_1)=d_A(u_0)P_A(u_1)+u_0\,d_A(P_A(u_1)) +\lambda d_A(u_0)d_A(P_A(u_1))= d(u_0)\ot u_1+u_0u_1+\lambda d(u_0)u_1.$$

\begin{theorem}$($\cite{GGZ,GK3}$)$
Let $Y$ be a set with a set map $d_0\colon Y\to \bfk[Y]$ and let $(\bfk[Y],d)$ be the
commutative differential algebra of weight $\lambda$ in
Theorem~\mref{thm:diff}.\mref{it:diffs}. The triple $(\sha
(\bfk[Y]), d_{\bfk[Y]},
    P_{\bfk[Y]})$, together with $j_{\bfk[Y]} \colon \bfk [Y] \to \sha(\bfk[Y])$, is the free commutative differential Rota-Baxter
    $\bfk$-algebra of weight $\lambda$ on the differential algebra $(\bfk[Y],d)$.
\mlabel{thm:freediffrb}
\end{theorem}

Apply Theorem~\mref{thm:freediffrb} to $Y:=\Delta X$ and $d_0$ as in Theorem~\mref{thm:diff}.\mref{it:commfreediff}. From Eq.~\ref{eq:rbm}, the set
\begin{equation}
\calb(\Delta X):=\left\{ u_0\ot \cdots \ot u_k\,\big|\, u_i\in
\cm(\Delta X), 0\leq i\leq k, k\geq 0\right\} \mlabel{eq:drbm}
\end{equation}
is a $\bfk$-basis of the free commutative differential Rota-Baxter
algebra $\sha(\Delta X)$ on the free differential algebra $(\bfk\{X\},d_X)$. We call this basis the set of {\bf differential
Rota-Baxter (DRB) monomials} on $X$.
Similarly, for $n\geq 1$, apply Theorem~\ref{thm:freediffrb} to  $Y:=\Delta_n X$ and $d_0$ as in Eq.~(\ref{eq:diffn}) of Theorem~\mref{thm:diff}.\ref{it:diffordn}. Then $\calb(\Delta_n X)$ is a
basis of $\sha(\Delta_n X)$ and is called the set of {\bf DRB
monomials of order $n$} on $X$.


\subsection{Free commutative operated algebras}
The construction of the free commutative operated algebra on a set $X$
that has the free commutative (differential) Rota-Baxter algebra as
a quotient is given in \mcite{GGZ}. The explicit construction $\sha(X)$ of the
free commutative Rota-Baxter algebra in Theorem~\mref{thm:shua} can
be realized on a submodule of the free commutative operated algebra
spanned by reduced words under a rewriting rule defined by the
Rota-Baxter axiom.

This construction is parallel to that of the free (noncommutative)
operated algebra on a set in~\mcite{BCQ,Gop,Gub,GSZ}.
See~\mcite{Qiu} for the non-unitary case.

\begin{defn}
{\rm Let $\Omega$ be a set. A {\bf commutative $\Omega$-operated monoid}
is a commutative monoid $G$ together with maps $\alpha_\omega\colon G\to
G,\omega\in \Omega$. A homomorphism between commutative $\Omega$-operated
monoids $(G,\{\alpha_\omega\}_\omega)$ and
$(H,\{\beta_\omega\}_\omega)$ is a monoid homomorphism $f\colon G\to H$
such that $f\circ \alpha_\omega=\beta_\omega\circ f$ for $\omega\in
\Omega$. }
\end{defn}

We similarly define the concept of a commutative $\Omega$-operated $\bfk$-algebra. The suffix $\Omega$ will be suppressed when the meaning of $\Omega$ is clear from the context.
We recall the construction of the free objects in the category of commutative
operated monoids \mcite{GGZ}.

Fix a set $Y$. Define monoids $\frakC_n:=\frakC_n(Y)$ for $n\geq
0$ by a recursion. First denote $\frakC_0:=\cm(Y)$.
For each $\omega\in \Omega$, let $\lc
\cm(Y)\rc_\omega:=\{\lc u\rc_\omega\,|\, u\in \cm(Y)\}$ be a set in bijection with $\cm(Y)$. We require that all the sets $\cm(Y)$ and $\lc \cm(Y)\rc_\omega, \omega \in \Omega$ are disjoint from one another. We write the notation $\sqcup$ for the disjoint union. Then define
$$\frakC_1:= \cm(Y\sqcup (\sqcup_{\omega\in \Omega} \lc \cm(Y)\rc_\omega))=\cm(Y\sqcup (\sqcup_{\omega\in \Omega} \lc \frakC_0\rc_\omega)).$$
Note that elements in $\lc \cm(Y)\rc_\omega$ are only symbols
indexed by elements in $\cm(Y)$. For example, $\lc 1\rc_\omega$ is
not the identity, but a new symbol. The inclusion $Y\hookrightarrow
Y\sqcup (\sqcup_{\omega\in \Omega}\lc\frakC_0\rc_\omega)$ induces a
monomorphism $i_{0}\colon \frakC_0=\cm(Y)\hookrightarrow
\frakC_1=\cm(Y\sqcup (\sqcup_{\omega\in \Omega}\, \lc\frakC_0\rc_\omega))$ of
free commutative monoids through which we identify $\frakC_0$ with
its image in $\frakC_1$. Inductively assume that $\frakC_{n-1}$ have
been defined for $n\geq 2$ and that the injection
$$i_{n-2}\colon \frakC_{n-2} \hookrightarrow \frakC_{n-1}$$
has been obtained. Then define
\begin{equation}
 \frakC_n:=\cm(Y\sqcup (\sqcup_{\omega\in \Omega} \lc\frakC_{n-1}\rc_\omega) ).
 \mlabel{eq:frakm}
 \end{equation}
Also the injection $i_{n-2}$ gives an injection
$$  \lc\frakC_{n-2}\rc_\omega \hookrightarrow
    \lc \frakC_{n-1} \rc_\omega, \ \omega\in \Omega.$$
Thus by the freeness of $\frakC_{n-1}=\cm(Y\sqcup (\sqcup_{\omega\in \Omega}
\lc\frakC_{n-2}\rc_\omega))$ as a free commutative monoid, we obtain
\begin{eqnarray*}
i_{n-1}\colon \frakC_{n-1} &=& \cm(Y\sqcup (\sqcup_{\omega\in \Omega}
\lc\frakC_{n-2}\rc_\omega))\hookrightarrow
    \cm(Y\sqcup (\sqcup_{\omega\in \Omega} \lc \frakC_{n-1}\rc_\omega)) =\frakC_{n}.
\end{eqnarray*}
Finally, define the commutative monoid
$$ \frakC(Y):=\bigcup_{n\geq 0}\frakC_n=\dirlim \frakC_n.$$
Elements in $\cmapm{Y}$ are called {\bf (commutative) $\Omega$-bracketed monomials} in $Y$.
Defining
\begin{equation}
\lc\ \rc_\omega\colon \cmapm{Y}\to \cmapm{Y}, \quad u\mapsto \lc u\rc_\omega, \
\omega\in \Omega, \mlabel{eq:mapp}
\end{equation}
then $(\cmapm{Y}, \{\lc\ \rc_\omega\}_{\omega\in \Omega})$ is a commutative operated
monoid and its linear span $(\bfk\cmapm{Y}, \{\lc\ \rc_\omega\}_{\omega\in \Omega})$ is a
commutative (unitary) $\Omega$-operated $\bfk$-algebra with its multiplication extended from $\frakC(Y)$ by linearity.

\begin{prop}(\mcite{GGZ}) Let $\Omega$ be a set.
\begin{enumerate}
\item
Let $j_Y\colon Y \hookrightarrow \cmapm{Y}$ be the natural embedding. Then the triple $(\cmapm{Y},\{\lc\ \rc_\omega\}_{\omega\in \Omega}, j_Y)$ is the free
commutative operated monoid on $Y$. More precisely, for any
commutative operated monoid $G$ and set map $f\colon Y\to G$, there is a
unique extension of $f$ to a homomorphism $\free{f}\colon\cmapm{Y}\to G$
of operated monoids.
\item
Let $j_Y\colon Y \hookrightarrow \bfk\cmapm{Y}$ be the natural embedding. Then the triple $(\bfk\cmapm{Y},\{\lc\ \rc_\omega\}_{\omega\in \Omega}, j_Y)$ is the
free commutative operated unitary $\bfk$-algebra on $Y$. More
precisely, for any commutative $\bfk$-algebra $R$ and set map
$f\colon Y\to R$, there is a unique extension of $f$ to a homomorphism
$\free{f}\colon\bfk\cmapm{Y}\to R$ of operated $\bfk$-algebras.
\end{enumerate}
\mlabel{pp:freetm}
\end{prop}

By the universal property of $\bfk\cmapm{Y}$, the
following conclusion from general principles of universal
algebra is obtained~\mcite{BN,Co}.

\begin{prop}(\mcite{GGZ})
Let $Y$ be a set with $d_0\colon Y\to \bfk[Y]$. Let $\Omega=\{d, P\}$ and write
$d(u):=\lc u\rc_d, P(u):=\lc u\rc_P$\,. Let $I_{DRB}=I_{DRB,Y}$ be
the operated ideal of $\bfk\cmapm{Y}$ generated by the set
$$\left\{\left . \begin{array}{l}
d(r)-d_0(r), \\
d(uv)-d(u)v-ud(v)-\lambda d(u)d(v),\\
 P(u)P(v)-P(uP(v)) -P(P(u) v) -\lambda P(uv),\\
 (d\circ P)(u)-u \end{array} \,\right|\, r\in Y, u, v\in \frak\cm(Y)\right\}.$$
Then the quotient operated algebra $\bfk\cmapm{Y}/I_{DRB}$, with operations induced by $d$ and $P$ (which we again denote by $d$ and $P$, respectively), is the free commutative
differential Rota-Baxter algebra on the differential algebra $(\bfk[Y],d)$ in Theorem~\mref{thm:diff}.\mref{it:diffs}. \mlabel{pp:freerb}
\end{prop}

Combining Proposition~\mref{pp:freerb} with Theorem~\mref{thm:freediffrb},
we have

\begin{prop}$($\mcite{GGZ}$)$
Let $Y$ be a set with $d_0\colon Y\to \bfk[Y]$. The natural embedding
$$\eta\colon \sha(\bfk[Y])\hookrightarrow \bfk\,\cmapm{Y}, \quad
u_0\ot u_1\ot \cdots \ot u_k \mapsto u_0P(u_1P( \cdots P(u_k) \cdots
)), \ u_i\in C(Y), 0\leq i\leq k, k\geq 0, $$ composed with the
quotient map $\rho:=\rho_Y\colon \bfk\,\frakC(Y) \to
\bfk\,\cmapm{Y}/I_{DRB}$ gives a linear bijection (in fact, an
isomorphism of differential Rota-Baxter algebras)
$$ \theta:=\theta_Y\colon \sha(\bfk[Y])\to \bfk\,\cmapm{Y}/I_{DRB}.$$
\mlabel{pp:comp}
\end{prop}

Because of the bijectivity of $\theta$, we can identify the basis $\calb(Y)$ of
$\sha(\bfk[Y])$ in Eq.~(\mref{eq:rbm}) with its image $\eta(\calb(Y))$ in
$\bfk\frakC(Y)$:
\begin{equation}
x_0\ot x_1 \ot \cdots \ot x_k  \leftrightarrow x_0P (x_1 P(\cdots
P(x_k)\cdots )),\ x_i\in \cm(Y), 0\leq i\leq k, k\geq 0.
\mlabel{eq:rbid}
\end{equation}

Define the {\bf reduction map}
\begin{equation}
\red:=\red_Y:=\theta^{-1}\circ \rho\colon \bfk\,\cmapm{Y} \to
\bfk\cmapm{Y}/I_{DRB,Y} \to \sha(\bfk[Y])\cong \eta(\sha(\bfk[Y])). \mlabel{eq:red}
\end{equation}
It reduces any $(d,P)$-bracketed monomial on $Y$ to a linear combination of DRB monomials on $Y$. For example, if
$u, v\in \cm(Y)$, then
$$\red(\lc u\rc_P \lc v\rc_P)=\red(P(u)P(v))=1\ot u\ot v + 1\ot v\ot u +\lambda \ot uv\leftrightarrow P( uP( v)) +P( vP( u)) +\lambda P(uv).$$


\section{Free commutative integro-differential algebras by initialization}
In this section, we summarize the construction of free commutative integro-differential algebras by initialization \mcite{GRR}.

\subsection{Regular differential algebras}
\mlabel{ss:reg-diffalg}
The construction applies to a large class of differential algebras called regular differential algebras. So we begin with the concept and examples of regular differential algebras.

\subsubsection{Quasi-antiderivatives and regularity}

\begin{defn}
Let $(A,\D)$ be a differential algebra of weight $\lambda$ with derivation $d$. A linear map $Q\colon A\to A$ is called a {\bf quasi-antiderivative} if
$d\circ Q \circ d=d$ and $Q \circ d \circ Q =Q $, with the additional condition that $\ker Q $ is a nonunitary $\bfk$-subalgebra of $A$ when $\lambda \neq 0$.
A differential algebra whose derivation has a quasi-antiderivative is called {\bf regular}.
\mlabel{de:reg}
\end{defn}
Given a regular differential ${\bf k}$-algebra $(A,\D)$ and a
fixed quasi-antiderivative $Q$ for $d$, we define the
following operators. Let
$$E= id_A-Q \circ d, S=d \circ
Q, J=id_A-E=Q \circ d, T=id_A - S.$$
We also
define $A_J$ to be the ${\bf k}$-submodule $A_J= im Q$,
and $A_T$ to be the ${\bf k}$-subalgebra $\ker Q$.

Regularity is equivalent to the existence of certain projectors, namely idempotent linear maps to a subspace.

\begin{prop}(\mcite{GRR})
Let $(A,\D)$ be a regular differential algebra. If $A$ is regular and $Q$ a quasi-antiderivation for $\D$, then the corresponding $\s: = d
\circ
  Q \colon A \to A$ is a projector onto~$\im\,d$ and~$\E := \id_A -
  Q  \circ d\colon A \to A$ is a projector onto~$\ker{d}$.
Conversely, if there are projectors $\s\colon A\to A$ onto $\im\, d$ and $\E\colon A\to
A$ onto $\ker d$, then there is a unique quasi-antiderivative $Q $ of $d$
such that $\im\, Q =\ker\,\E$ and $\ker Q =\ker \s$ and $(A,\D)$ is regular. \label{pp:reg}
\end{prop}

To prove the converse, from the given projectors $S$ and $E$, we have $A=\ker S \oplus \im\,d$ and $A=\ker\,d \oplus \ker\,E$. Thus the restriction of $d$ to $\ker\,E$ is a bijection onto $\im\,d$. Then there is unique map $Q \colon A\to A$ whose restriction to $\im\,d$ is the inverse of the above bijection and whose kernel is $\ker S$.

\subsubsection{Differential polynomial algebras}
Let $Y$ be a set with a well-ordering $\leq_Y$. Define the {\bf length-lexicographic order} $\leq_{Y,\lex}^*$ on the free monoid $M(Y)$ by
\begin{equation}
u<_{Y,\lex}^* v \Leftrightarrow \left\{\begin{array}{l} \ell< m, \\
\text{or } \ell=m \text{ and } \exists 1\leq i_0\leq \ell \text{ such that } u_i=v_i \text{ for } 1\leq i<i_0 \text{ and } u_{i_0}<v_{i_0}, \end{array}\right.
\mlabel{eq:lex}
\end{equation}
where $u=u_1\cdots u_\ell$ and $v=v_1\cdots v_m$ with $u_i,v_j\in Y, 1\leq i\leq \ell, 1\leq j\leq m.$.
It is well-known~\mcite{BN} that $\leq_{Y,\lex}^*$ is again a well-ordering.
An element $1\neq u$ of the free commutative monoid $\cm(Y)$ can be uniquely expressed as
\begin{equation}
u = u_0^{j_0} \cdots u_k^{j_k}, \text{ where } u_0,\cdots,u_k\in
Y, j_0,\cdots, j_k \in \mathbb{Z}_{\geq 1} \text{ and } u_0
> \cdots > u_k. \mlabel{eq30}
\end{equation}
This expression is called the {\bf standard form} of $u$. If $k= -1$, we take $u=1$ by convention.

Let $X$ be a well-ordered set and let $Y=\Delta X$ (resp. $\Delta_n X$). Let $n\geq 0$ be
given. For $x_0^{(i_0)}, x_1^{(i_1)}\in Y$ with $x_0, x_1\in X$, define
\begin{equation}
x_0^{(i_0)} \leq x_1^{(i_1)} \left(\text{resp. } x_0^{(i_0)}\leq_n
x_1^{(i_1)}\right) \Leftrightarrow (x_0,-i_0) \leq (x_1, -i_1) \quad
\text{ lexicographically}. \mlabel{eq:difford}
\end{equation}
For example $x^{(2)} < x^{(1)}< x$. Also, $x_1<x_2$ implies
$x_1^{(i_1)} < x_2^{(i_2)}$ for all $i_1, i_2\geq 0$.

\begin{defn}
Let $u\in \cm(\Delta X)$ with standard form in Eq.~(\mref{eq30}):
\begin{align}
u=u_0^{j_0}\cdots u_k^{j_k}, \text{ where } u_0,\cdots, u_k\in
\Delta X, u_0 > \cdots > u_k \text{ and } j_0,\cdots,j_k\in
\mathbb{Z}_{\geq 1}. \notag
\end{align}
Call $u$ {\bf functional} if either $u=1$ or $u_k\in X$ or $j_k >1$.
\end{defn}

\begin{prop}(\mcite{GGZ,GRR})
  \label{pp:diffpoly}
Let $\lambda \in {\bf k}$ and let $X$ be a set. Let $A=({\bf k}\{X\}, d_\lambda)$ be the free commutative
differential algebra of weight $\lambda$ on $X$ as defined in Theorem~\mref{thm:diff}.\mref{it:commfreediff}. Then there are direct sums~$A =  A_T \oplus \im\,{\D}$ and~$A = A_J \oplus \ker{\D}$, where
\begin{equation}
\cala_T=\cala_{T,n}=\{u\in C(\Delta_n X)\,|\, u \text{ is functional}\}, \quad  A_T := \bfk \cala_T,
\mlabel{eq:polyf}
\end{equation}
and $A_J$ is the submodule generated by all monomials $1\ne u \in C(\Delta X)$.
Thus $d$ admits a quasi-antiderivative~$\q$. Therefore, $(\bfk\{X\},d_\lambda)$ is regular.
\end{prop}
Since the product of two functional
monomials is again functional, $A_T$ is in fact a $\bfk$-subalgebra of $A$.

As noted in the remark after Proposition~\mref{pp:reg}, the quasi-antiderivative~$\q$ is defined as follows. From the direct sums, the derivation $D$ restricts to a bijection $D\colon A_J\to \im\,{\D}$. Define $\q\colon \im\,{\D}\to A_J$ to be the inverse map and then extend $\q$ to $A$ by taking $A_T$ to be the kernel of $\q$.

\subsubsection{Rational functions}

We show that the algebra of rational functions with derivation of any weight is regular.

Let $A=\CC(x)$. For given $\lambda\in \CC$, let
\begin{equation*}
\D_\lambda \colon A\to A, f(x)\mapsto \left\{\begin{array}{ll}
\frac{f(x+\lambda)-f(x)}{\lambda}, &\lambda\neq 0, \\
f'(x), & \lambda =0,
\end{array} \right .
\end{equation*}
be the $\lambda$-derivation.
Denote

$$
\calr\colon = \left\{\begin{array}{ll}
\left\{ \sum\limits_{i=1}^k \sum\limits_{j=1}^{n_i} \frac{\gamma_{ij}}{(x-\alpha_i)^j}\, \Big|\, \alpha_i\in \CC \text{ distinct }, \gamma_{ij}\in \CC\right\}, & \lambda =0, \\
& \\
\left\{ \sum\limits_{i=1}^k \sum\limits_{j=1}^{n_i}
    \frac{\gamma_{ij}}{(x-\alpha_{ij})^i}\,\Big|\,\alpha_{ij}\in \CC \text{ distinct for any given } i, \gamma_{ij}\in \CC \text{ nonzero}\right\}, & \lambda \neq 0. \end{array} \right.
$$
Then denote
$$ \CC(x)_J \colon= x\CC[x] + \calr$$
and
$$ \CC(x)_T\colon=\left \{\begin{array}{ll}
\left \{\sum\limits_{i=1}^k \frac{\gamma_i}{x-\alpha_i}\,\Big|\, \alpha_i\in \CC \text{ distinct }, \gamma_i\in \CC\right\}, & \lambda =0, \\
& \\
\left \{ \sum\limits_{i=1}^k \sum\limits_{j=1}^{n_i}
    \frac{\gamma_{ij}}{(x-\alpha_{ij})^i}\in \calr\,\Big| \mathrm{re}(\alpha_{ij})\in [0,|\mathrm{re} (\lambda)|) \right\}, & \lambda\neq 0, \end{array} \right .
$$
where $\mathrm{re}(z)$ is the real part of $z\in \CC$.
It is proved in~\mcite{GRR} that
\begin{equation}
\CC(x)=\im\, \D_\lambda \oplus \CC(x)_T.
\mlabel{eq:fracdecomp1}
\end{equation}
Further, $\CC(x)_T$ is a nonunitary subalgebra of $\CC(x)$.
We also have
\begin{equation*}
 \CC(x)=\ker \D_\lambda \oplus \CC(x)_J.
 \end{equation*}
Then by Proposition~\mref{pp:reg}, $\D_\lambda$ is regular.

\subsection{Construction of $\efid{A}$}
\mlabel{ss:freec}
We now give the construction of the free commutative integro-differential algebra $\efid{A}$ on a regular differential algebra $(A, \D)$ with a fixed quasi-antiderivative $\q$.

With the notations set up after Definition~\ref{de:reg}, we give now an explicit construction of $\efid{A}$ via free commutative Rota-Baxter algebras and tensor products. First let
\begin{equation*}
  \sha_\ct (A):= \bigoplus_{k\geq 0} A\ot A_\ct^{\otimes k} = A\oplus (A\ot A_\ct) \oplus (A\ot A_\ct^{\ot 2}) +\cdots
\end{equation*}
be the $\bfk$-submodule of $\sha(A)$ in Eq.~(\mref{eq:freerb}). Then
$\sha_\ct(A)$ is the tensor product $A\ot \sha^+(A_\ct)$
where $\sha^+(A_\ct):=\bigoplus\limits_{n\geq 0} A_\ct^{\ot n}$ is
the mixable shuffle algebra~\mcite{Gub,GK1,Ho} on the non-unitary $\bfk$-algebra $A_\ct$.

Next, let $K:=\ker\, d \supset \bfk$ and let
\begin{equation*}
  A_\ee:=\{\ee(a)\,|\, a\in A\}
\end{equation*}
denote a replica of the $K$-algebra $A$, endowed with the zero derivation and the $K$-algebra structure map
$$K\to A_\ee,\quad  c\mapsto \ee(c), c\in K.$$
We will use the $K$-algebra isomorphism
$$\ee\colon A \rightarrow A_\ee,\quad a\mapsto \ee(a), a\in A.$$
Let
\begin{equation}
  \efid{A}:=A_\ee\ot_K \sha_\ct(A)=A_\ee\ot_K A\ot \sha^+(A_\ct)
  \mlabel{eq:efid}
\end{equation}
denote the tensor product differential algebra of $A_\ee$ and $\sha_\ct(A)$, namely the tensor product algebra where the derivation $d_A$ is defined by the Leibniz rule. To define the linear operator $\pp_A$ on $\efid{A}$, we first require that $\pp_A$
be $A_\ee$-linear. Then we just need to define $\pp_A(\fraka)$
for a pure tensor $\fraka$ in $A\ot \sha^+(A_\ct)$. For this purpose we apply induction on the length $n$ of $\fraka$.

When $n=1$, we have $\fraka = a \in A$. By definition of $T$ we have
$a = \D(\q(a)) + \T(a)$ with $\T(a) \in A_\ct$. Then we define
\begin{equation}
  \pp_A(a) := \q(a) -\ee(\q(a)) + 1 \ot \T(a).
\mlabel{eq:pa}
\end{equation}
Assume that $\pp_A(\fraka)$ has been defined for $\fraka$ of length $n\geq 1$ and consider the case when $\fraka$ has length $n+1$. Then
$\fraka=a\ot \lbar{\fraka}$ where $a\in A, \lbar{\fraka}\in
A_\ct^{\ot n}$ and we define
\begin{equation}
  \pp_A(a\ot \lbar{\fraka}):= \q(a) \ot \lbar{\fraka} -
  \pp_A(\q(a) \lbar{\fraka})
  -\lambda \, \pp_A(\D(\q(a)) \, \lbar{\fraka})+ 1 \ot \T(a) \ot
  \lbar{\fraka},
  \mlabel{eq:Pu}
\end{equation}
where the first and last terms are already in $A\ot
\sha^+(A_\ct)$ while the middle terms are in $\efid{A}$ by the
induction hypothesis.

\begin{theorem}(\mcite{GRR})
  Let~$(A, \D)$ be a regular differential algebra of
  weight~$\lambda$ with a fixed quasi-anti\-derivative~$\q$. Then the triple $(\efid{A},
  \D_A, \pp_A)$, with the natural embedding
  $$i_A\colon  A \hookrightarrow \efid{A}=A_\vep\ot_K A\ot \sha^+(A_\ct)$$
  onto the
  second tensor factor, is the free commutative integro-differential
  algebra of weight $\lambda$ generated by $A$.
  \mlabel{thm:intdiffa}
\end{theorem}

\section{Free commutative integro-differential algebras by Gr\"obner-Shirshov bases}

In this section, we give a construction of free commutative integro-differential algebras by the method of Gr\"obner-Shirshov bases. The main result Theorem~\mref{thm:gsb} can be read independently of the rest of the section, which is meant to give some details of the method.

The method of Gr\"obner bases or Gr\"obner-Shirshov bases originated
from the work of Buchberger~\mcite{Bu} (for commutative polynomial
algebras, 1965), Hironaka~\mcite{Hi} (for infinite series algebras, 1964),
Shirshov~\mcite{Sh} (for Lie algebras, 1962) and Zhukov~\mcite{Zh} (reduction in nonassociative algebra, 1950). It has since become a fundamental method in commutative algebra, algebraic geometry and computational algebra, and has been extended to many other algebraic structures, notably associative algebras~\mcite{Be,Bo}. In recent
years, the method of Gr\"obner-Shirshov bases has been applied to a
large number of algebraic structures to study problems on normal
forms, word problems, rewriting systems, embedding theorems,
extensions, growth functions and Hilbert series.
See~\mcite{BCC,BCL,BLSZ} for further details.

The method of Gr\"obner bases or Gr\"obner-Shirshov bases is very useful in constructing free objects in various categories, including the
alternative constructions of free Rota-Baxter algebras and free
differential Rota-Baxter algebras~\mcite{BCD,BCQ}. The basic idea is
to prove a composition-diamond lemma that achieves a rewriting
procedure for reducing any element to a certain ``standard form". Then
the set of elements in standard form is a basis of the free object.

In the recent paper\mcite{GGZ}, this method is applied to construct the free
commutative integro-differential algebra as the quotient of a free commutative differential Rota-Baxter algebra modulo the integration by parts formula in Eq.~(\mref{eq:ibp1}). In order to do so, the authors first establish a Composition-Diamond Lemma for the free
commutative differential Rota-Baxter algebra constructed
in~\mcite{GK3}. Then they prove that the ideal generated by the
defining relation of integro-differential algebras in
Eq.~(\mref{eq:ibp1}) has a Gr\"obner-Shirshov basis, thereby
identifying a basis of the free commutative integro-differential
algebra as a canonical subset of a known basis of a free
commutative differential Rota-Baxter algebra.

\subsection{Weakly monomial order} \mlabel{subsec:mon} In this
subsection, we will define a weak form of the monomial order on pieces of the set of differential Rota-Baxter monomials filtered by the order of differentiation. It
will be sufficient to establish the composition-diamond lemma
for integro-differential algebras.

For a set $X$, recall that $\Delta X \colon=\{x^{(k)}
\mid x\in X, k\geq 0\}$ and $\Delta_n X \colon=\{x^{(k)} \mid x\in X, 0\leq k\leq n\}$
for $n\geq 0$. Then the family $\{\cm(\Delta_n X)\}_{n\geq 0}$ defines an increasing filtration on $\cm(\Delta X)$ and hence by Eq.~(\mref{eq:rbm}), induces a filtration $\{\calb(\Delta_n X)\}_{n\geq 0}$ of the set $\calb(\Delta X)$ of DRB monomials by DRb monomials of order $n$.
Elements of $\calb(\Delta_n X)$ are called
{\bf DRB monomials of order $n$}.

In Definition~\mref{de:drb} below and what follows, the DRB ($\star$-DRB) monomials are elements in the basis $\calb(\Delta X)$ (resp. $\calb(\Delta X^\star)$) of $\sha(\bfk\{X\})$ (resp. $\sha(\bfk\{X^\star\})$), which are identified via Eq.~(\ref{eq:rbid}) as $(d,P)$-bracketed monomials $\eta(\calb(\Delta X))\subseteq \bfk\frakC(\Delta X)$ (resp. $\eta(\calb(\Delta X^{\star}))\subseteq \bfk \frakC(\Delta X^\star)$).

\begin{defn}
{\rm Let $X$ be a set, $\star$ a symbol not in $X$ and $\Delta_n
X^\star := \Delta_n (X\cup \{\star\})$.
\begin{enumerate}
\item
By a {\bf
$\star$-DRB monomial on $\Delta_n X$},
we mean any expression in $\calb(\Delta_n X^\star)$ with exactly one
occurrence of $\star$. The set of all $\star$-DRB monomials on $\Delta_n X$ is denoted by
$\calb^\star(\Delta_n X)$.
\item
For $q\in \calb^\star(\Delta_n X)$ and
$u\in \calb(\Delta_n X)$, we define
$$q|_u := q|_{\star \mapsto u}$$
to be the bracketed monomial in $\cmapm{\Delta_n X}$ obtained by replacing the letter $\star$ in $q$ by
$u$, and call $q|_u$ a {\bf $u$-monomial on $\Delta_n X$}.
\item
Further, for $s=\sum_i c_i u_i \in \bfk
\calb(\Delta_n X)$, where $c_i\in \bfk$, $u_i\in \calb(\Delta_n X)$
and $q\in \calb^\star(\Delta_n X)$, we define
$$q|_s := \sum_i c_i q|_{u_i},$$
which is in $\bfk\,\cmapm{\Delta_n X}$.
\end{enumerate}
\mlabel{de:drb}
}
\end{defn}
We note that a $\star$-DRB monomial $q$ is a DRB monomial in $\Delta_n X^\star$ while its substitution $q|_u$ might not be a DRB monomials. For example, for $q=P(x_1)\star\in \eta(\calb(\Delta_n X^\star))$ and $u=P(x_2)\in \calb(\Delta_n X)$ where $x_1, x_2\in X$, the $u$-monomial $q|_u=P(x_1)P(x_2)$ is no longer in $\eta(\calb(\Delta_n X))$.

\begin{defn}
If $q = p|_{d^\ell(\star)}$ for some $p\in \calb^\star(\Delta_n
X)$ and $\ell\in \mathbb{Z}_{\geq 1}$, then we call $q$ a {\bf
type I $\star$-DRB monomial}. Let $\calb_{\rm I}^\star(\Delta_n X)$ denote the set of type I $\star$-DRB monomials on
$\Delta_n X$ and call
$$\calb_{\rm II}^\star(\Delta_n X) :=
\calb^\star(\Delta_n X) \setminus \calb_{\rm I}^\star(\Delta_n X)$$
the set of {\bf type II $\star$-DRB monomials}.
\end{defn}

For example, $d(\star) P(x)\in \calb_{\rm I}^\star(\Delta_n X)$ and $\star  P(x) \in \calb_{\rm II}^\star(\Delta_n X)$.

\begin{defn}{\rm
Let $X$ be a set, $\star_1$, $\star_2$ two distinct symbols not in
$X$ and $\Delta_n X^{\star_1, \star_2} := \Delta_n (X\cup
\{\star_1,\star_2\})$. We define a {\bf
$(\star_1,\star_2)$-DRB monomial on
$\Delta_n X$ } to be an expression in $\calb(\Delta_n
X^{\star_1,\star_2})$ with exactly one occurrence of $\star_1$ and
exactly one occurrence of $\star_2$. The set of all $(\star_1,
\star_2)$-DRB monomials on $\Delta_n X$
is denoted by $\calb^{\star_1, \star_2}(\Delta_n X)$. For $q\in
\calb^{\star_1, \star_2}(\Delta_n X)$ and $u_1, u_2\in \bfk
\calb(\Delta_n X)$, we define
$$q|_{u_1,u_2} := q|_{\star_1 \mapsto u_1, \star_2 \mapsto u_2}$$
to be the bracketed monomial obtained by
replacing the letter $\star_1$ (resp. $\star_2$) in $q$ by $u_1$
(resp. $u_2$) and call it a {\bf $(u_1,u_2)$-bracketed monomial on $\Delta_n X$ }. }
\end{defn}

A $(u_1,u_2)$-DRB monomial on $\Delta_n
X$ can also be recursively defined by
\begin{equation*} \label{eq12}
q|_{u_1,u_2} := (q^{\star_1}|_{u_1})|_{u_2},
\end{equation*}
where $q^{\star_1}$ is $q$ when $q$ is regarded as a
$\star_1$-DRB monomial on the set
$\Delta_n X^{\star_2}$. Then $q^{\star_1}|_{u_1}$ is in
$\calb^{\star_2}(\Delta_n X)$. Similarly, we have
\begin{equation*} \label{eq13}
q|_{u_1,u_2} := (q^{\star_2}|_{u_2})|_{u_1}.
\end{equation*}

Let $X$ be a well-ordered set. Let $n\geq 0$ be
given. We extend the well-ordering $\leq_n$ on $\cm(\Delta_n X)$ defined in Eq.~(\mref{eq:difford}) to $\calb(\Delta_n X)$. Note that $$\calb(\Delta_n X)=
\{u_0\ot u_1\ot \cdots \ot u_k\,|\, u_i\in \cm(\Delta_n X), 1\leq
i\leq k, k\geq 0\}=\sqcup_{k\geq 1} \cm(\Delta_nX)^{\ot k} $$ can be
identified with the free semigroup on the set $\cm(\Delta_nX)$. Thus
the well-ordering $\leq_n$ on $\cm(\Delta_nX)$ extends to a well-ordering $\leq_{n,\lex}^*$~\mcite{BN} on $\calb(\Delta_n X)$ which we will still denote by $\leq_n$
for simplicity. More precisely, for any $u=u_0\otimes \cdots \otimes
u_k \in \cm(\Delta_n X)^{\ot (k+1)}$ and $v = v_0 \otimes \cdots
\otimes v_\ell \in \cm(\Delta_n X)^{\ot (\ell+1)}$, define
\begin{equation}
u \leq_{n} v \text { if }(k, u_0, \cdots,u_k) \leq (\ell, v_0,
\cdots,v_\ell) \text{ lexicographically}. \mlabel{eq7}
\end{equation}

\begin{defn}
Let $\leq_n$ be the well-ordering on $\calb(\Delta_nX)$ defined in
Eq.~(\mref{eq7}). Let $q\in \calb^\star(\Delta_n X)$ and $s\in
\mathbf{k}\calb(\Delta_n X)$.
\begin{enumerate}
\item
For any $0\neq f\in \bfk \calb(\Delta_n X)$, let $\lbar{f}$ denote
the {\bf leading term} of $f$: $f = c \overline{f} + \sum_{i} c_iu_i$,
where $0\neq  c, c_i\in \bfk$, $u_i\in \calb(\Delta_n X)$, $u_i<
\overline{f}$. We call $f$ {\bf monic} if $c=1$.
\item
Let
$$\overline{q |_s} := \overline{\red(q|_{s}}),$$
where $\red\colon \bfk\frakC(\Delta_nX)\to \sha(\Delta_nX)=\eta(\bfk
\calb(\Delta_n X))$ is the reduction map in Eq.~(\mref{eq:red}).
\item
The element $q|_s\in \bfk\,\cmapm{\Delta_nX}$ is called {\bf normal}
if $q|_{\lbar{s}}$ is in $\calb(\Delta_n X)$. In other words, if
$\red(q|_{\lbar{s}}) = q|_{\lbar{s}}$.
\end{enumerate}
\mlabel{normaldef}
\end{defn}

\begin{remark}
By definition, $q|_s$ is normal if and only if
$q|_{\lbar{s}}$ is normal if and only if the $\lbar{s}$-DRB monomial
$q|_{\lbar{s}}$ is already a DRB monomial, that is, no further
reduction in $\sha(\Delta_n X)$ is possible.
\end{remark}

Here are some examples of abnormal $s$-DRB monomials.
\begin{exam}
\begin{enumerate}
\item $q=\star P(y)$ and $\bar{s}=P(x)$, giving $q|_{\bar{s}}=P(x)P(y)$ which is reduced to $P(xP(y))+P(P(x)y)+\lambda P(xy)$ in $\eta(\sha(\Delta_n X))$;
\item $q=d(\star)$ and $\bar{s}=P(x)$, giving $q|_{\bar{s}}=d(P(x))$ which is reduced to $x$ in $\eta(\sha(\Delta_n X))$;
\item $q=d(\star)$ and $\bar{s}=x^2$, giving $q|_{\bar{s}}=d(x^2)$ which is reduced to $2xx^{(1)}+\lambda (x^{(1)})^2$ in $\eta(\sha(\Delta_n X))$;
\item $q=d^n(\star)$ and $\bar{s}=d(x)$, giving $q|_{\bar{s}}=d^{n+1}(s)$ which is reduced to $0$ in $\eta(\sha(\Delta_n X))$.
\end{enumerate}
\end{exam}

\begin{defn}\label{defweakmonomial}
A {\bf weakly monomial order} on $\calb(\Delta_n X)$ is a well-ordering
$\leq$ satisfying the following condition:
\begin{quote}
For $u, v\in \calb(\Delta_n X)$, if $u \leq v$, then $\overline{q|_u} \leq \overline{q|_v}$  if $q \in
\calb^{\star}_{\rm II}(\Delta_n X)$, or if $q \in
\calb^{\star}_{\rm I}(\Delta_n X)$ and $q|_v$  is normal.
\end{quote}
\end{defn}

\begin{prop}$($\mcite{GGZ}$)$
The order $\leq_n$ defined in Eq.~(\ref{eq7}) is a weakly monomial order on $\calb(\Delta_n X)$. \mlabel{weakmonomial}
\end{prop}

\subsection{Composition-Diamond lemma} \mlabel{sec:cd} In this
section, we shall establish the composition-diamond lemma for the free commutative differential Rota-Baxter algebra
$\sha(\bfk[\Delta_n X])$ of order $n$.

\begin{defn}
\begin{enumerate}
\item Let $u,w \in \calb(\Delta_n X)$. We call $u$ a {\bf
subword} of $w$ if $w$ is in the operated ideal of $\frakC(\Delta_n
X)$ generated by $u$. In terms of $\star$-words, $u$ is a subword of
$w$ if there is a $q\in \calb^\star(\Delta_n X)$ such that $w=q|_u$.
\item
Let $u_1$ and $u_2$ be two subwords of $w$. Then $u_1$ and $u_2$ are
called {\bf separated} if $u_1\in \cm(\Delta_n X)$, $u_2\in
\calb(\Delta_n X)$ and there is a $q\in \calb^{\star_1,
\star_2}(\Delta_n X)$ such that $w=q|_{u_1,u_2}$.
\item
For any $u\in \calb(\Delta_n X)$, $u$ can be expressed as $u =u_{1}
\cdots u_{k}$, where $u_1, \cdots, u_{k-1} \in \Delta_n X$ and
$u_{k} \in \Delta_n X\cup P(\calb(\Delta_n X))$. The expression is unique up to permutations of those factors in $\Delta_n X$. The integer $k$ is called
the {\bf breadth} of $u$ and is denoted by $\mathrm{bre}(u)$.
\item
Let $f,g\in \calb(\Delta_n X)$. A pair $(u,v)$ with $u\in
\calb(\Delta_n X)$ and $v\in \cm(\Delta_n X)$ is called an {\bf
intersection pair} for $(f,g)$ if the differential
Rota-Baxter monomial $w:= fu$ equals $vg$ and satisfies
$\mathrm{bre}(w)<\mathrm{bre}(f) + \mathrm{bre}(g)$. Then we call
$f$ and $g$ to be {\bf overlapping}. Note that if $f$ and $g$ are
overlapping, then $f\in \cm(\Delta_n X)$.
\end{enumerate}
 \mlabel{lemma:overlape}
\end{defn}
For example, let $w=xyxy$ with $x, y\in X$ and $u_1=xy$ be the subword of $w$ on the left and $u_1=xy$ be the subword of $w$ on the right. Then $u_1$ and $u_2$ are separated. Let $g$ be the subword $yx$ of $w$. Then $u_1$ and $g$ are overlapping. A systematic discussion on relative locations (separated, overlapping and inclusion) of two subwords can be found in~\mcite{GGSZ,GZ2}.

There are three kinds of compositions.

\begin{defn}
{\rm
Let $\leq_n$ be the weakly monomial order on $\calb(\Delta_n X)$ defined in Eq.~(\ref{eq7}), and let
$f,g\in \mathbf{k}\calb(\Delta_n X)$ be monic with respect to $\leq_n$ such that $f\neq g$.

\begin{enumerate}
\item If $\overline{f}\in C( \Delta_n X) P(\calb(\Delta_n X))$, then define a {\bf composition of multiplication} to be $fu$ where $u\in C(\Delta_n X) P(\calb(\Delta_n X))$.
\item If there is an intersection pair $(u,v)$ for $(\overline{f},
\overline{g})$, then we define
$$(f,g)_w :=(f,g)^{u,v}_{w} := fu-vg$$ and call it an {\bf intersection
composition} of $f$ and $g$.
\item If there exists a $q\in \calb^\star(\Delta_n X)$ such that
$w:=\overline{f}=q|_{\overline{g}}$, then we define $(f,g)_w
:=(f,g)^{q}_{w} := f-q|_g$ and call it an {\bf inclusion
composition} of $f$ and $g$ with respect to $q$. Note that if this
is the case, then $q|_g$ is normal.
\end{enumerate}
}
\end{defn}
In the last two cases, $(f,g)_w$ is called the {\bf ambiguity} of the composition. For example, let $$f=P(d(u)P(d(v)P(r)))-uP(d(v)P(r))+P(ud(v)P(r))+\lambda P(d(u)d(v)P(r))$$
and
$$g=P(d(v)P(r))-vP(r)+P(vr)+\lambda P(d(v)r)$$
with the first terms being the leading terms. Then
we have $\bar{f}=q|_{\bar{g}}$ where $q\colon = P(d(u)\star)$. Hence we get an inclusion composition of $f$ and $g$ with the ambiguity
\begin{eqnarray*}
(f,g)_w^q&=& -uP(d(v)P(r))+P(ud(v)P(r))+\lambda P(d(u)d(v)P(r))\\
&&-\left ( -P(d(u)vP(r))+P(d(u)P(vr))+\lambda P(d(u)P(d(v)r))\right).
\end{eqnarray*}

\begin{defn}
{\rm Let $\leq_n$ be the weakly monomial order on $\calb(\Delta_n X)$ defined in Eq.~(\mref{eq7}),
$S\subseteq \mathbf{k}\calb(\Delta_n X)$ be a set of monic
differential Rota-Baxter polynomials and $w\in
\calb(\Delta_n X)$.

\begin{enumerate}
\item A composition of multiplication $fu$ is called {\bf trivial mod
$[S]$} if $$fu = \sum_i c_i q_i|_{s_i}, $$ where $c_i\in \bfk$,
$q_i\in \calb^\star(\Delta_n X)$, $s_i\in S$, $q_i|_{s_i}$ is normal
and $q_i|_{\overline{s_i}} \leq_n \overline{fu}$. If this is the case,
we write $$fu \equiv 0  \text{ mod }  [S].$$
\item For $u,v\in \mathbf{k}\calb(\Delta_n X)$ and $w\in \calb(\Delta_n X)$, we say $u$ and $v$
are {\bf congruent modulo $(S,w)$} and denote this by $$u \equiv v \text
{ mod } (S,w)$$ if $u-v = \sum_{i}c_iq_i|_{s_i}$, where $c_i\in
\mathbf{k}$, $q_i\in \calb^\star(\Delta_n X)$, $s_i\in S$,
$q_i|_{s_i}$ is normal and $\overline{q_i|_{s_i}} <_n w$.
\item For $f,g\in \mathbf{k}\calb(\Delta_n X)$ and suitable $u,v$ or $q$
that give an intersection composition $(f,g)^{u,v}_{w}$ or an inclusion composition $(f,g)^{q}_{w}$, the composition is called
{\bf trivial modulo $(S,w)$} if $$(f,g)^{u,v}_{w} \text{ or }
(f,g)^{q}_{w}\equiv 0 \text{ mod } (S,w).$$
\item The set $S\subseteq \mathbf{k}\calb(\Delta_n X)$ is a {\bf Gr\"{o}bner-Shirshov basis} if all
compositions of multiplication are trivial mod $[S]$, and, for $f, g\in S$, all
intersection compositions $(f,g)^{u,v}_{w}$ and all inclusion
compositions $(f,g)^{q}_{w}$ are trivial modulo $(S,w)$.
\end{enumerate}
}
\end{defn}

\begin{theorem}(\mcite{GGZ} Composition-Diamond Lemma) Let $\leq_n$ be the weakly monomial order on $\calb(\Delta_n X)$
defined in Eq.~(\ref{eq7}), $S_n$ a set of monic DRB polynomials in $\bfk \calb(\Delta_n
X)$ with $d(S_n) \subseteq S_n$, and $\mathrm{Id}(S_n)$ the Rota-Baxter ideal of
$\mathbf{k}\calb(\Delta_n X)$ generated by $S_n$. Then with respect to $\leq_n$, the following
conditions are equivalent:
\begin{enumerate}
\item $S_n$ is a Gr\"{o}bner-Shirshov basis in
$\mathbf{k}\calb(\Delta_n X)$.
\mlabel{it:cda}
\item If $0\neq f\in \mathrm{Id}(S_n)$, then $\overline{f} =
q|_{\overline{s}}$ for some $q\in \calb^\star(\Delta_n X)$, $s\in S_n$
and $q|_s$ is normal.
\mlabel{it:cdb}
\item
The set $\mathrm{Irr}(S_n):= \calb(\Delta_n X) \setminus
\{ q|_{\overline{s}} \mid q\in \calb^\star(\Delta_n X), s\in S_n ,
q|_{s}\text{ is normal} \}$ is a $\mathbf{k}$-basis of
$\mathbf{k}\calb(\Delta_n X)/\mathrm{Id}(S_n)$.
In other words, $\mathbf{k} \mathrm{Irr}(S_n) \oplus \mathrm{Id}(S_n) = \mathbf{k}
\calb(\Delta_n X)$.
\mlabel{it:cdc}
\end{enumerate} \mlabel{thm:CD lemma}
\end{theorem}

\subsection{Free commutative integro-differential algebras by Gr\"obner-Shirshov bases} \mlabel{subsec:gs}

In this subsection we begin with a finite set $X$ and prove that the
relation ideal of the free commutative differential Rota-Baxter
algebra on $X$ of order $n\geq 1$, defining the
corresponding commutative integro-differential algebra of order $n$
possesses a Gr\"obner-Shirshov basis. Then we apply the
Composition-Diamond Lemma in Theorem~\mref{thm:CD lemma} to
construct a canonical basis for the free commutative integro-differential
algebra of order $n$. As $n$ approaches infinity, we obtain a
canonical basis of the free commutative integro-differential algebra
on the finite set $X$. Finally for any well-ordered set $X$, by
showing that the canonical basis of the free commutative
integro-differential algebra on each finite subset of $X$ is
compatible with the inclusion of the subset in $X$, we obtain a
canonical basis of the free commutative integro-differential algebra
on $X$.

\begin{theorem}(\mcite{GGZ})
Let
\begin{equation}
 S_n:= \left\{P(d(u) P(v))- uP(v)+ P(uv) + \lambda P(d(u) v)\,\big|\, u, v\in  \sha(\bfk[\Delta_n X])\right\} \label{eq2}
\end{equation}
be the set of generators corresponding to the integration by parts axiom Eq.~(\mref{eq:ibp1}).
Let $\leq_n$ be the monomial order defined in Eq.~(\mref{eq7}).
\begin{enumerate}
\item
With respect to $\leq_n$, $S_n$ is a Gr\"{o}bner-Shirshov basis in $\sha(\bfk[\Delta_n X])$. Hence $\mathrm{Irr}(S_n)$ in Theorem~\mref{thm:CD lemma} is a linear basis of $\sha(\bfk[\Delta_n X])/\Id(S_n)$.
\item
Let $A_T=\bfk\{X\}_T$ be as defined in Eq.~(\mref{eq:polyf}), $A_n = \bfk [\Delta_n X]$, $A_{n,T} =
A_{n}\cap A_T$. Let $I_{{\rm ID},n}$ be the differential Rota-Baxter ideal
of $\sha(A_n)$ generated by $S_n$. Then we have the isomorphism of modules
$$\sha(A_n) / I_{ID,n} \cong A_n \oplus \left(\bigoplus_{k\geq 0} A_n\otimes A_{n-1,T}^{\otimes k} \otimes A_n \right).$$
\end{enumerate}
\mlabel{thm:PBWBase}
\end{theorem}

Let
\begin{equation}
S  := \left\{P(d(u) P(v))- uP(v)+ P(uv) + \lambda P(d(u) v) \mid
u, v\in \sha(\Delta X) \right\}. \mlabel{eq:gsid}
\end{equation}
be the set of generators corresponding to the integration by parts axiom Eq.~(\mref{eq:ibp1}).

\begin{theorem}(\mcite{GGZ})
Let $X$ be a nonempty well-ordered set, $A_T=\bfk\{X\}_T$,
$\sha(\bfk\{X\})=\sha(\Delta X)$ the free commutative
differential Rota-Baxter algebra on $X$ and $I_{ID}$ the ideal
of $\sha(\bfk\{X\})$ generated by $S$ defined in
Eq.~(\mref{eq:gsid}). Then the composition
$$\sha(A)_T:=A \oplus \left(\bigoplus_{k\geq 0} A\otimes A_{T}^{\otimes k} \otimes A \right) \hookrightarrow  \sha(A) \to \sha(A) / I_{ID}$$
of the inclusion and the quotient map is an isomorphism of $\bfk$-modules.
 \mlabel{thm:gsb}
\end{theorem}

It would be interesting to compare the two constructions of free commutative integro-differential algebras in Theorem~\mref{thm:intdiffa} and Theorem~\mref{thm:gsb}. The advantages of the first construction is that it applies to a large class of differential algebras and that the product in the free algebra is clearly defined. The advantage of the second construction is that the construction comes from a subset of the free commutative differential Rota-Baxter algebra from which the free integro-differential algebra is obtained modulo an ideal. It is useful to have both of the two constructions available in order to study different aspects of free commutative integro-differential algebras. Further study in this direction is being pursued in another work. The construction of free noncommutative integro-differential algebras is also under investigation.

\smallskip

\noindent
{\bf Acknowledgements}:
This work is supported by the National Natural Science Foundation of China (Grant No. 11201201 and 11371178), Fundamental Research Funds for the Central Universities (Grant No. lzujbky-2013-8), the Natural Science Foundation of Gansu Province (Grant No. 1308RJZA112) and the National Science Foundation of US (Grant No. DMS~1001855). The authors thank the referees and editors for helpful comments.

\end{document}